\documentclass[11pt]{article}
\input epsf.tex
\usepackage{amsmath,amssymb,latexsym,times}
\catcode`\@=11
\renewcommand\subsection{\@startsection{subsection}{2}{\z@}%
                                     {-3.25ex\@plus -1ex \@minus -.2ex}%
                                     {-0.01 mm}
                                     {\normalfont\large\bfseries}}

\catcode`\@=11
\renewcommand\subsubsection{\@startsection{subsubsection}{2}{\z@}%
                                     {-3.25ex\@plus -1ex \@minus -.2ex}%
                                     {-0.01 mm}
                                     {\normalfont\bfseries}}

\newtheorem{example}{Example}
\newtheorem{theorem}{Theorem}
\newtheorem{definition}{Definition}
\newtheorem{proposition}{Proposition}
\newtheorem{lemma}{Lemma}

\def\resp{{\em resp.$\ $}}

\def\proof{\medskip\noindent {\it Proof --- \ }}

\def\cqfd{\hfill $\Box$ \bigskip}
\def\adots{\mathinner{\mkern2mu\raise1pt\hbox{.}
\mkern3mu\raise4pt\hbox{.}\mkern1mu\raise7pt\hbox{.}}}
\def\<{\langle\,}
\def\>{\,\rangle}

\def\hom{{\rm Hom}}

\def\ext{{\rm Ext}}
\def\End{{\rm End}}

\def\l{\lambda}
\def\a{\alpha}
\def\aa{\mathbf a}
\def\de{\delta}
\def\b{\beta}
\def\bb{\mathfrak{b}}

\def\N{{\mathbb N}}

\def\C{{\mathbb C}}

\def\F{{\cal F}}

\def\P{{\mathbb P}}

\def\g{\mathfrak g}
\def\h{\mathfrak h}

\def\Sl{\mathfrak{sl}}

\def\bb{\mathfrak b}

\def\n{{\mathfrak n}}

\def\G{{\cal G}}

\def\L{\Lambda}

\def\<{\langle}
\def\>{\rangle}

\def\E{{\cal E}}

\def\deg{{\rm deg}}

\def\M{{\cal M}}

\def\le{\leqslant}
\def\ge{\geqslant}

\def\De{\Delta}

\def\i{\mathbf i}

\def\lra{\longrightarrow}
\def\lla{\longleftarrow}

\def\1{\mathbf 1}

\def\ra{\rightarrow} 
\def\dimv{\mathbf{dim}}

\newcommand{\bsm}{\begin{smallmatrix}}
\newcommand{\esm}{\end{smallmatrix}}
\def\f{\mathfrak f}
\def\0{\mathbf{0}}
\def\socle{{\rm socle}}

\def\id{{\rm id}}
\def\nm{\mathfrak{n}_{-}}
\def\oM{\widetilde{\cal M}}
\def\pg{\mathfrak{p}}
\def\Pg{\mathfrak{P}}
\def\qg{\mathfrak{q}}
\def\Qg{\mathfrak{Q}}

\newdimen\Squaresize \Squaresize=14pt
\newdimen\Thickness \Thickness=0.5pt
 
\def\Square#1{\hbox{\vrule width \Thickness
   \vbox to \Squaresize{\hrule height \Thickness\vss
      \hbox to \Squaresize{\hss#1\hss}
   \vss\hrule height\Thickness}
\unskip\vrule width \Thickness}
\kern-\Thickness}
 
\def\Vsquare#1{\vbox{\Square{$#1$}}\kern-\Thickness}

\title{\bf Verma modules and preprojective algebras}

\author{Christof {\sc Geiss}
\thanks{C. Geiss acknowledges support from DGAPA grant IN101402-3.}
, Bernard {\sc Leclerc} 
\thanks{B. Leclerc is grateful to the GDR 2432 and the GDR 2249
for their support.}\ \ 
and Jan {\sc Schr\"oer}
\thanks{J. Schr\"oer was supported by a research 
fellowship from the DFG (Deutsche Forschungsgemeinschaft).}}

\date{}

\begin{document}
\maketitle
\begin{abstract}
We give a geometric construction of the Verma modules of a
symmetric Kac-Moody Lie algebra $\g$ in terms of constructible
functions on the varieties of nilpotent finite-dimensional
modules of the corresponding preprojective algebra $\L$.
\end{abstract}

\section{Introduction}
Let $\g$ be the symmetric Kac-Moody Lie algebra associated 
to a finite unoriented graph $\Gamma$ without loop.
Let $\nm$ denote a maximal nilpotent subalgebra of $\g$.
In \cite[\S 12]{Lu91}, Lusztig has given a geometric construction
of $U(\nm)$ in terms of certain
Lagrangian varieties. 
These varieties can be interpreted as module varieties
for the preprojective algebra $\L$ attached to the graph 
$\Gamma$ by Gelfand and Ponomarev \cite{GP}.
In Lusztig's construction, $U(\nm)$ gets identified with
an algebra $(\M,*)$ of constructible functions on these
varieties, where $*$ is a convolution product 
inspired by Ringel's multiplication for Hall algebras.

Later, Nakajima gave a similar construction of the highest
weight irreducible integrable $\g$-modules $L(\l)$ in terms
of some new Lagrangian varieties which differ from Lusztig's
ones by the introduction of some extra vector spaces $W_k$ for
each vertex $k$ of $\Gamma$, and by considering only stable
points instead of the whole variety \cite[\S 10]{Na}.

The aim of this paper is to extend Lusztig's original construction
and to endow $\M$ with the structure of a Verma module $M(\l)$.
 
To do this we first give a variant of the geometrical construction
of the integrable $\g$-modules $L(\l)$, using
functions on some natural open subvarieties 
of Lusztig's varieties instead of functions
on Nakajima's varieties (Theorem~\ref{thI}).
These varieties have a simple description in terms of
the preprojective algebra $\L$ and of certain injective
$\L$-modules $q_\l$. 

Having realized the integrable modules $L(\l)$ as quotients
of $\M$, it is possible, using the comultiplication of $U(\nm)$,
to construct geometrically the raising operators $E_i^\l\in\End(\M)$ 
which make $\M$ into the Verma module $M(\l)$ (Theorem~\ref{conjV}).
Note that we manage in this way to realize Verma modules with
arbitrary highest weight (not necessarily dominant).

Finally, we dualize this setting and give a geometric construction of 
the dual Verma module $M(\l)^*$ in terms of the delta functions
$\de_x \in \M^*$ attached to the finite-dimensional nilpotent
$\L$-modules $x$ 
(Theorem~\ref{dual}).

\section{Verma modules}
\label{sect1}
\subsection{}
Let $\g$ be the symmetric Kac-Moody Lie algebra associated 
with a finite unoriented graph $\Gamma$ without loop.
The set of vertices of the graph is denoted by~$I$.
The (generalized) Cartan matrix of $\g$ is $A=(a_{ij})_{i,j\in I}$,
where $a_{ii}=2$ and, for $i\not = j$, $-a_{ij}$ is the number
of edges between $i$ and $j$.

\subsection{}
Let $\g = \n\oplus \h\oplus \nm$ be a Cartan decomposition of
$\g$, where $\h$ is a Cartan subalgebra and $(\n,\nm)$ a 
pair of opposite maximal
nilpotent subalgebras.
Let $\bb=\n\oplus\h$.
The Chevalley generators of $\n$ (\resp $\nm$) are denoted
by $e_i\ (i\in I)$ (\resp $f_i$) and we set $h_i=[e_i,f_i]$. 

\subsection{}
Let $\alpha_i$ denote the simple root
of $\g$ associated with $i\in I$.
Let $(-\,;\,-)$ be a symmetric bilinear form on $\h^*$ such
that $(\a_i\,;\,\a_j)=a_{ij}$.  
The lattice of integral weights in $\h^*$ is denoted by $P$,
and the sublattice spanned by the simple roots is denoted
by $Q$.
We put
\[
P_+=\{\l\in P \mid (\l\,;\,\a_i) \ge 0, \ (i\in I)\}, 
\qquad Q_+=Q\cap P_+.
\]

\subsection{}
Let $\l\in P$ and let $M(\lambda)$ be the Verma module 
with highest weight $\l$.
This is the induced $\g$-module defined by
$M(\l) = U(\g)\otimes_{U(\bb)} \C\,u_\l$,
where $u_\l$ is a basis of the one-dimensional representation of 
$\bb$ given by
\[
h\,u_\l = \l(h)\,u_\l,\quad n\,u_\l = 0,\qquad (h\in \h,\ n\in\n).
\]
As a $P$-graded vector space $M(\l)\cong U(\nm)$ 
(up to a degree shift by $\l$). 
$M(\l)$ has a unique simple quotient denoted by $L(\l)$,
which is integrable if and only if $\l\in P_+$.
In this case, the kernel of the $\g$-homomorphism
$M(\l) \ra L(\l)$ is the $\g$-module $I(\l)$
generated by the vectors 
\[
f_i^{(\l\,;\,\a_i)+1}\otimes u_\l,\qquad (i\in I).
\]

\section{Constructible functions}

\subsection{}
Let $X$ be an algebraic variety over $\C$ endowed with its Zariski topology.
A map $f$ from $X$ to a vector space $V$ is said to
be constructible if its image $f(X)$ is finite, and for each $v\in
f(X)$ the preimage $f^{-1}(v)$ is a constructible subset of $X$.

\subsection{}
By $\chi(A)$ we denote the Euler characteristic of a constructible 
subset $A$ of $X$.
For a constructible map $f : X \ra V$ one defines
\[
\int_{x\in X} f(x) = \sum_{v\in V}  \chi(f^{-1}(v))\,v \in V.
\]
More generally, for a constructible subset $A$ of $X$ we write
\[
\int_{x\in A} f(x) = \sum_{v\in V} \chi(f^{-1}(v)\cap A)\,v.
\]

\section{Preprojective algebras}

\subsection{}
Let $\Lambda$ be the preprojective algebra associated to the 
graph $\Gamma$ (see for example \cite{Ri,GLS}).
This is an associative $\C$-algebra, which is
finite-dimensional if and only if 
$\Gamma$ is a graph of type $A, D, E$.
Let $s_i$ denote the simple one-dimensional $\Lambda$-module
associated with $i\in I$, and let $p_i$ be its projective cover
and $q_i$ its injective hull.
Again, $p_i$ and $q_i$ are finite-dimensional if and only if 
$\Gamma$ is a graph of type $A, D, E$.

\subsection{}
A finite-dimensional $\Lambda$-module $x$ is nilpotent if and
only if it has a composition series with all factors of the
form $s_i\ (i\in I)$. 
We will identify the dimension vector of $x$ with an element
$\b\in Q_+$ by setting $\dimv(s_i) = \a_i$.

\subsection{}\label{embed}
Let $q$ be an injective $\L$-module of the form
\[
q = \bigoplus_{i\in I} q_i^{\oplus a_i}
\]
for some nonnegative integers $a_i\ (i\in I)$.

\begin{lemma}
Let $x$ be a finite-dimensional $\L$-module isomorphic to
a submodule of $q$. If $f_1 : x \ra q$ and $f_2 : x \ra q$
are two monomorphisms, then 
there exists an automorphism $g : q \ra q$ such that
$f_2 = gf_1$.
\end{lemma}
\proof
Indeed, $q$ is the injective hull of its socle
$b=\bigoplus_{i\in I}\,s_i^{\oplus a_i}$.
Let $c_j\ (j=1,2)$ be a complement of
$f_j(\socle(x))$ in $b$.
Then $c_1\cong c_2$ and the maps
\[
h_j := f_j \oplus \id : \quad x \oplus c_j \ra q, \qquad (j=1,2)
\]
are injective hulls.
The result then follows from the unicity of the injective hull.
\cqfd

Hence, up to isomorphism, there is a unique way
to embed $x$ into $q$.

\subsection{}
Let $\M$ be the algebra of constructible functions on the varieties
of finite-dimensional nilpotent $\Lambda$-modules defined by 
Lusztig \cite{Lu00} to give a geometric realization of $U(\nm)$.
We recall its definition. 

For $\b=\sum_{i\in I}b_i\a_i\in Q_+$, 
let $\L_\b$ denote the variety of nilpotent $\L$-modules 
with dimension vector $\b$.
Recall that $\L_\b$ is endowed with an action of the algebraic
group $G_\b = \prod_{i\in I}GL_{b_i}(\C)$, so that two points
of $\L_\b$ are isomorphic as $\L$-modules if and only if
they belong to the same $G_\b$-orbit.
Let $\oM_\b$ denote the vector space of constructible
functions from $\L_\b$ to $\C$ which are constant on $G_\b$-orbits.
Let
\[
\oM=\bigoplus_{\b\in Q_+} \oM_\b.
\]
One defines a multiplication $*$ on $\oM$ as follows.
For $f\in \oM_\b$, $g\in\oM_\gamma$ and $x\in \L_{\b+\gamma}$,
we have 
\begin{equation}\label{sta}
(f*g)(x) = \int_{U} f(x') g(x''),
\end{equation}
where the integral 
is over the variety of $x$-stable subspaces
$U$ of $x$ of dimension $\gamma$, $x''$ is the $\L$-submodule 
of $x$ obtained by restriction to $U$ and $x'=x/x''$. 
In the sequel in order to simplify notation, we will not
distinguish between the subspace $U$ 
and the submodule $x''$ of $x$ carried by~$U$.
Thus we shall rather write 
\begin{equation}\label{star}
(f*g)(x) = \int_{x''} f(x/x'') g(x''),
\end{equation}
where the integral 
is over the variety of submodules
$x''$ of $x$ of dimension $\gamma$.

For $i\in I$, the variety $\L_{\a_i}$ is reduced to a single 
point : the simple module $s_i$.
Denote by $\1_i$ the function mapping this point to $1$.
Let $\G(i,x)$ denote the  variety of all submodules
$y$ of $x$ such that $x/y \cong s_i$.
Then by (\ref{star}) we have
\begin{equation}\label{stari}
(\1_i*g)(x) = \int_{y\in\G(i,x)} g(y).
\end{equation}

Let $\M$ denote the subalgebra of $\oM$ generated by the
functions $\1_i\ (i\in I)$.
By Lusztig~\cite{Lu00}, $(\M,*)$ is isomorphic to $U(\nm)$
by mapping $\1_i$ to the Chevalley generator $f_i$.

\subsection{}
In the identification of $U(\nm)$ with $\M$, formula 
(\ref{stari}) represents the left multiplication by $f_i$.
In order to endow $\M$ with the structure of a Verma module
we need to introduce the following important definition.
For $\nu \in P_+$, let 
\[
q_\nu = \bigoplus_{i\in I} q_i^{\oplus (\nu\,;\,\a_i)}.
\]
Lusztig has shown \cite[\S 2.1]{Lu002} that Nakajima's Lagrangian varieties
for the geometric realization of $L(\nu)$ are isomorphic to the
Grassmann varieties of $\L$-submodules of $q_\nu$ with a
given dimension vector.

Let $x$ be a finite-dimensional nilpotent $\L$-module isomorphic to
a submodule of the injective module $q_\nu$.
Let us fix an embedding $F : x \ra q_\nu$ and identify $x$
with a submodule of $q_\nu$ via $F$. 
 
\begin{definition}\label{def}
For $i\in I$ let 
$\G(x,\nu,i)$ be the variety of submodules $y$ of $q_\nu$ containing $x$ and 
such that $y/x$ is isomorphic to $s_i$.
\end{definition}

This is a projective variety which, by \ref{embed},
depends only (up to isomorphism) on $i$, $\nu$ and the 
isoclass of $x$.

\section{Geometric realization of integrable irreducible $\g$-modules}
\label{gri}

\subsection{}
For $\l\in P_+$ and $\b\in Q_+$,
let $\L_{\b}^{\l}$ denote the variety of nilpotent $\L$-modules of dimension
vector $\b$ which are isomorphic to a submodule of $q_\l$.
Equivalently $\L_{\b}^{\l}$ consists of the nilpotent modules of dimension
vector $\b$ whose socle contains $s_i$ with multiplicity at most $(\l\,;\,\a_i)$
$(i\in I)$.
This variety has been considered by Lusztig \cite[\S 1.5]{Lu03}.
In particular it is known that $\L_{\b}^{\l}$ is an open subset of
$\L_{\b}$, and that the number of its irreducible components
is equal to the dimension of the $(\l-\b)$-weight space of $L(\l)$.

\subsection{}
Define $\oM_\b^\l$ to be the vector space of constructible functions
on $\L_{\b}^{\l}$ which are constant on $G_\b$-orbits.
Let $\M_\b^\l$ denote the subspace of $\oM_\b^\l$
obtained by restricting elements of $\M_\b$ to 
$\L_{\b}^{\l}$.
Put $\oM^\l=\bigoplus_\b\oM_\b^\l$ and
$\M^\l=\bigoplus_\b\M_\b^\l$.
For $i\in I$ define endomorphisms $E_i, F_i, H_i$ of $\oM^\l$ 
as follows:
\begin{eqnarray}
(E_if)(x) &=& \int_{y\in\G(x,\l,i)} f(y),
\qquad (f\in \oM_\b^\l,\ x\in \L_{\b-\a_i}^{\l}),\label{actE}\\[3mm]
(F_if)(x) &=& \int_{y\in\G(i,x)} f(y),
\quad\qquad (f\in \oM_\b^\l,\ x\in \L_{\b+\a_i}^{\l}),\label{actF} \\[3mm]
(H_if)(x) &=& (\l-\b;\a_i)\,f(x),
\qquad (f\in \oM_\b^\l,\ x\in \L_{\b}^{\l}).\label{actH}
\end{eqnarray}

\begin{theorem}\label{thI}
The endomorphisms $E_i, F_i, H_i$ of $\oM^\l$ leave stable the
subspace $\M^\l$. Denote again by $E_i, F_i, H_i$ the induced
endomorphisms of $\M^\l$.
Then the assignments $e_i\mapsto E_i$, $f_i\mapsto F_i$, $h_i\mapsto H_i$,
give a representation of $\g$ on $\M^\l$ isomorphic to the
irreducible representation $L(\l)$.
\end{theorem}

\subsection{}
The proof of Theorem~\ref{thI} will involve a series of lemmas.

\subsubsection{}
For $\i=(i_1,\ldots ,i_r)\in I^r$ and 
    $\aa=(a_1,\ldots ,a_r)\in \N^r$, define the 
variety $\G(x,\l,(\i,\aa))$ of flags of $\L$-modules
\[
\f = (x=y_0 \subset y_1 \subset \cdots \subset y_r \subset q_\l)
\]
with $y_k/y_{k-1} \cong s_{i_k}^{\oplus a_k}\ (1\le k\le r)$.
As in Definition~\ref{def}, this is a projective variety 
depending (up to isomorphism) only on $(\i,\aa)$, $\l$ and the isoclass of $x$
and not on the choice of a specific embedding of
$x$ into $q_\l$.
\begin{lemma}\label{prod}
Let $f\in \oM_\b^\l$ and $x\in
\L_{\b-a_1\a_{i_1}-\cdots-a_r\a_{i_r}}^{\l}$.
Put $E_i^{(a)}=(1/a!)E_i^a$.
We have
\[
(E_{i_r}^{(a_r)}\cdots E_{i_1}^{(a_1)}f)(x) = 
\int_{\f\in\G(x,\l,(\i,\aa))} f(y_r).
\]
\end{lemma}
The proof is standard and will be omitted.

\subsubsection{}
By \cite[12.11]{Lu91} the endomorphisms $F_i$ satisfy the
Serre relations
\[
\sum_{p=0}^{1-a_{ij}} (-1)^p\, F_j^{(p)}\, F_i\, F_j^{(1-a_{ij}-p)}=0
\]
for every $i\not = j$.
A similar argument shows that
\begin{lemma}\label{serre}
The endomorphisms $E_i$ satisfy the
Serre relations
\[
\sum_{p=0}^{1-a_{ij}} (-1)^p\, E_j^{(p)}\, E_i\, E_j^{(1-a_{ij}-p)}=0
\]
for every $i\not = j$.
\end{lemma}
\proof
Let $f\in\oM_\b^\l$ and $x\in \L^\l_{\b- \a_i - (1-a_{ij})\a_j}$.
By Lemma~\ref{prod},
\[
(E_j^{(p)}\, E_i\, E_j^{(1-a_{ij}-p)}f)(x) =
\int_\f f(y_3)
\]
the integral being taken on the variety of flags
\[
\f = (x \subset y_1 \subset y_2 \subset y_3 \subset q_\l)
\]
with $y_1/x \cong s_j^{\oplus 1-a_{ij}-p}$, $y_2/y_1 \cong s_i$
and $y_3/y_2 \cong s_j^{\oplus p}$.
This integral can be rewritten as
\[
\int_{y_3} f(y_3)\,\chi(\F[y_3;p])
\]
where the integral is now over all submodules $y_3$ of
$q_\l$ of dimension $\b$ containing $x$  
and $\F[y_3;p]$ is the variety of flags $\f$
as above with fixed last step $y_3$.
Now, by moding out the submodule $x$ at each step of the flag,
we are reduced to the same situation as in \cite[12.11]{Lu91},
and the same argument allows to show that 
\[
\sum_{p=0}^{1-a_{ij}} \chi(\F[y_3;p]) =0,
\]
which proves the Lemma.
\cqfd

\subsubsection{}
Let $x\in\L^\l_\b$. 
Let $\varepsilon_i(x)$ denote the multiplicity of $s_i$
in the head of $x$.
Let $\varphi_i(x)$ denote the multiplicity of $s_i$
in the socle of $q_\l/x$.
\begin{lemma}\label{keylem}
Let $i,j\in I$ (not necessarily distinct). 
Let $y$ be a submodule of $q_\l$ containing $x$ and 
such that $y/x \cong s_j$.
Then
\[
\varphi_i(y)-\varepsilon_i(y) =
\varphi_i(x)-\varepsilon_i(x) - a_{ij}.
\]
\end{lemma}
\proof
We have short exact sequences
\begin{eqnarray}
0 &\to& x\ \ \to\ \ \ q_\l\ \ \ \to\ \ q_\l/x\  \to\ \,0,
\label{eqeq1}\\
0 &\to& y\ \ \to\ \ \ q_\l\ \ \ \to\ \ q_\l/y\  \to\ \,0,
\label{eqeq1p}\\
0 &\to& x\  \ \to\ \ \ \ y\ \ \ \ \to\ \ \ s_j\ \ \ \,\to\ \ 0, \label{eqeq2}\\
0 &\to& s_j\ \to\ \ q_\l/x\ \to\ q_\l/y\ \to\ 0.  \label{eqeq3}
\end{eqnarray}
Clearly, $\varepsilon_i(x)=|\hom_\L(x,s_i)|$, the dimension of
$\hom_\L(x,s_i)$.
Similarly $\varepsilon_i(y)=|\hom_\L(y,s_i)|$,
$\varphi_i(x)=|\hom_\L(s_i,q_\l/x)|$,
$\varphi_i(y)=|\hom_\L(s_i,q_\l/y)|$.
Hence we have to show that
\begin{equation}\label{eqeq4}
|\hom_\L(x,s_i)|-|\hom_\L(y,s_i)| = 
|\hom_\L(s_i,q_\l/x)| - |\hom_\L(s_i,q_\l/y)| - a_{ij}.
\end{equation}
In our proof, we will use 
a property of preprojective algebras proved
in \cite[\S 1]{CB}, namely,
for any finite-dimensional
$\L$-modules $m$ and $n$ there holds 
\begin{equation}\label{eqCB}
|\ext_\L^1(m,n)|=|\ext_\L^1(n,m)|.
\end{equation}

(a)\ \ If $i=j$ then $a_{ij}=2$, $|\hom_\L(s_j,s_i)|=1$ and
$|\ext^1_\L(s_j,s_i)|=0$ since $\Gamma$ has no loops.
Applying $\hom_\L(-,s_i)$ to (\ref{eqeq2}) we get the exact sequence
\[
0 \to \hom_\L(s_j,s_i) \to \hom_\L(y,s_i) \to \hom_\L(x,s_i) \to 0,
\]
hence 
\[
|\hom_\L(x,s_i)|-|\hom_\L(y,s_i)| = -1.
\] 
Similarly applying $\hom_\L(s_i,-)$ to (\ref{eqeq3}) we get an exact sequence
\[
0 \to \hom_\L(s_i,s_j) \to \hom_\L(s_i,q_\l/x) \to \hom_\L(s_i,q_\l/y) \to 0,
\]
hence 
\[
|\hom_\L(s_i,q_\l/x)| - |\hom_\L(s_i,q_\l/y)|=1,
\] 
and (\ref{eqeq4}) follows.

(b)\ \ If $i\not = j$, we have $|\hom_\L(s_i,s_j)|=0$ and 
$|\ext^1_\L(s_i,s_j)|=|\ext^1_\L(s_j,s_i)|=-a_{ij}$.
Applying $\hom_\L(s_i,-)$ to (\ref{eqeq2}) we get an exact sequence
\[
0 \to \hom_\L(s_i,x) \to \hom_\L(s_i,y) \to 0,
\]
hence 
\begin{equation}\label{eqC}
|\hom_\L(s_i,x)|-|\hom_\L(s_i,y)| = 0.
\end{equation}
Moreover, by \cite[\S 1.1]{Bo},
$|\ext^2_\L(s_i,s_j)|=0$ because there are no relations from
$i$ to $j$ in the defining relations of $\L$.
(Note that the proof of this result in \cite{Bo} only requires 
that $I \subseteq J^2$
(here we use the notation of \cite{Bo}).
One does not need the additional assumption $J^n \subseteq I$ for some $n$. 
Compare also the discussion in \cite{BK}.)

Since $q_\l$ is injective
$|\ext^1_\L(s_i,q_\l)|=0$,
thus applying $\hom_\L(s_i,-)$ to  (\ref{eqeq1}) we get an exact 
sequence
\[
0\to \hom_\L(s_i,x) \to \hom_\L(s_i,q_\l)\to \hom_\L(s_i,q_\l/x) 
\to \ext^1_\L(s_i,x) \to 0,
\]
hence 
\begin{equation}\label{eqA}
|\hom_\L(s_i,x)|-|\hom_\L(s_i,q_\l)|+|\hom_\L(s_i,q_\l/x)|
-|\ext^1_\L(s_i,x)|=0.
\end{equation}
Similarly, applying $\hom_\L(s_i,-)$ to  (\ref{eqeq1p}) we get 
\begin{equation}\label{eqB}
|\hom_\L(s_i,y)|-|\hom_\L(s_i,q_\l)|+|\hom_\L(s_i,q_\l/y)|
-|\ext^1_\L(s_i,y)|=0.
\end{equation}
Subtracting (\ref{eqA}) from (\ref{eqB}) and taking into account
(\ref{eqCB}) and (\ref{eqC}) we obtain
\begin{equation}\label{eqD}
|\ext^1_\L(x,s_i)|-|\ext^1_\L(y,s_i)|=|\hom_\L(s_i,q_\l/x)|-|\hom_\L(s_i,q_\l/y)|.
\end{equation}
Now applying $\hom_\L(-,s_i)$ to  (\ref{eqeq2}) we get 
the long exact sequence
\[
0 \to \hom_\L(y,s_i) \to \hom_\L(x,s_i) \to  \ext^1_\L(s_j,s_i)
 \to  \ext^1_\L(y,s_i)  \to  \ext^1_\L(x,s_i) \to 0,
\]
hence
\[
|\hom_\L(y,s_i)|-|\hom_\L(x,s_i)|-a_{ij}
-|\ext^1_\L(y,s_i)|+|\ext^1_\L(x,s_i)|=0,
\]
thus, taking into account (\ref{eqD}), we have proved (\ref{eqeq4}).
\cqfd

\begin{lemma}\label{keylem2}
With the same notation we have 
\[
\varphi_i(x)-\varepsilon_i(x) = (\l-\b;\a_i).
\]
\end{lemma}
\proof
We use an induction on the height of $\b$.
If $\b=0$ then $x$ is the zero module and $\varepsilon_i(x) = 0$.
On the other hand $q_\l/x = q_\l$ and $\varphi_i(x)=(\l;\a_i)$
by definition of $q_\l$.
Now assume that the lemma holds for $x\in \L^\l_\b$ and let 
$y\in \L^\l_{\b+\a_j}$ be a submodule
of $q_\l$ containing $x$. 
Using Lemma~\ref{keylem} we get that 
\[
\varphi_i(y)-\varepsilon_i(y) = (\l-\b;\a_i)- a_{ij} =
(\l-\b-\a_j;\a_i),
\]
as required, and the lemma follows.
\cqfd

\begin{lemma}\label{lemH}
Let $f\in\oM^\l_\b$.
We have 
\[
(E_iF_j - F_jE_i)(f) = \de_{ij}(\l-\b ; \a_i)f.
\] 
\end{lemma}
\proof
Let $x\in\L^\l_{\b-\a_i+\a_j}$.
By definition of $E_i$ and $F_j$ we have
\[
(E_iF_jf)(x) = \int_{\pg\in\Pg} f(y)
\]
where $\Pg$ denotes the variety of pairs $\pg=(u,y)$
of submodules of $q_\l$ with $x \subset u$, $y \subset u$,
$u/x\cong s_i$ and $u/y\cong s_j$.
Similarly,
\[
(F_jE_if)(x) = \int_{\qg\in\Qg} f(y)
\]
where $\Qg$ denotes the variety of pairs $\qg=(v,y)$
of submodules of $q_\l$ with $v \subset x$, $v \subset y$,
$x/v\cong s_j$ and $y/v\cong s_i$.

Consider a submodule $y$ such that there exists in $\Pg$
(\resp in $\Qg$) at least one pair of the form $(u,y)$
(\resp $(v,y)$).
Clearly, the subspaces carrying the submodules $x$ and $y$
have the same dimension $d$ and their intersection has dimension
at least $d-1$. If this intersection has dimension exactly
$d-1$ then there is a unique pair $(u,y)$ 
(\resp  $(v,y)$), namely $(x+y,y)$ (\resp  $(x\cap y,y)$).
This means that
\[
\int_{\pg\in\Pg;\ y\not = x} f(y)
=
\int_{\qg\in\Qg;\ y\not = x} f(y).
\]
In particular, since when $i\not = j$ we cannot have $y=x$,
it follows that
\[
 (E_iF_j - F_jE_i)(f) = 0, \qquad (i\not = j).
\]
On the other hand if $i=j$ we have
\[
((E_iF_i - F_iE_i)(f))(x) = f(x) (\chi(\Pg')-\chi(\Qg'))
\]
where $\Pg'$ is the variety of submodules $u$ of $q_\l$
containing $x$ such that $u/x \cong s_i$,
and $\Qg'$ is the variety of submodules $v$ of $x$
such that $x/v \cong s_i$.
Clearly we have $\chi(\Qg')=\varepsilon_i(x)$
and $\chi(\Pg')=\varphi_i(x)$. 
The result then follows from
Lemma~\ref{keylem2}.
\cqfd

\subsubsection{}
The following relations for the endomorphisms
$E_i, F_i, H_i$ of $\oM^\l$ are easily checked
\[
[H_i,H_j]=0, \quad [H_i,E_j]=a_{ij}E_j,
\quad [H_i,F_j]=-a_{ij}F_j.
\]
The verification is left to the reader.
Hence, using Lemmas~\ref{serre} and \ref{lemH}, we have proved
that  the assignments $e_i\mapsto E_i$, $f_i\mapsto F_i$, $h_i\mapsto H_i$,
give a representation of $\g$ on $\oM^\l$.

\begin{lemma}
The endomorphisms
$E_i, F_i, H_i$ 
leave stable the subspace $\M^\l$.
\end{lemma}
\proof
It is obvious for $H_i$, and it follows from the definition
of $\M^\l$ for $F_i$.
It remains to prove that if $f\in\M^\l_\b$ then
$E_if\in\M^\l_{\b-\a_i}$.
We shall use induction on the height of $\b$.
We can assume that $f$ is of the form $F_jg$ for some
$g\in \M^\l_{\b-\a_j}$.
By induction we can also assume that 
$E_ig\in\M^\l_{\b-\a_i-\a_j}$.
We have
\[
E_if=E_iF_jg=F_jE_ig + \de_{ij}(\l- \b+\a_j;\a_i)g,
\]
and the right-hand side clearly belongs to $\M^\l_{\b-\a_i}$.
\cqfd

\begin{lemma}\label{lem7}
The representation of $\g$ carried by $\M^\l$ is isomorphic
to $L(\l)$.
\end{lemma}
\proof
For all $f\in\M_\b$ and all 
$x\in\L^\l_{\b+(a_i+1)\a_i}$ we have 
$f*\1_i^{*(a_i+1)}(x)=0$.
Indeed, by definition of $\L^\l$ the socle of $x$ contains
$s_i$ with multiplicity at most $a_i$.
Therefore the left ideal of $\M$ generated by the functions
$\1_i^{*(a_i+1)}$ is mapped to zero by the linear map 
$\M \ra \M^\l$ sending a function $f$ on $\L_\b$ to its 
restriction to $\L^\l_\b$. 
It follows that for all $\b$ the dimension of $\M^\l_\b$
is at most the dimension of the $(\l-\b)$-weight space
of $L(\l)$.

On the other hand, the function $\1_0$ mapping the zero $\L$-module
to $1$ is a highest weight vector of $\M^\l$ of weight $\l$.
Hence $\1_0\in\M^\l$ generates a quotient of the Verma module $M(\l)$,
and since $L(\l)$ is the smallest quotient of $M(\l)$ we must
have $\M^\l=L(\l)$.
\cqfd

This finishes the proof of Theorem~\ref{thI}.

\section{Geometric realization of Verma modules}\label{grv}

\subsection{}
Let $\b\in Q_+$ and $x \in \L_{\b-\a_i}$. 
Let $q = \bigoplus_{i\in I} q_i^{\oplus a_i}$
be the injective hull of $x$. 
For every $\nu\in P_+$ such that $(\nu;\a_i)\ge a_i$ the injective
module $q_\nu$ contains a submodule isomorphic to $x$.
Hence, for such a weight $\nu$ and for any $f\in\M_\b$, the integral 
\[
\int_{y\in\G(x,\nu,i)} f(y)
\]
is well-defined.
\begin{proposition}\label{conjf}
Let $\l\in P$ and choose $\nu\in P_+$ such that $(\nu;\a_i)\ge a_i$
for all $i\in I$.
The number
\begin{equation}\label{stabil}
\int_{y\in\G(x,\nu,i)} f(y)\ -\ (\nu - \l\,;\, \a_i)\, f(x\oplus s_i)
\end{equation}
does not depend on the choice of $\nu$.
Denote this number by $(E_i^\l f)(x)$. 
Then, the function 
\[
E_i^\l f : x \mapsto (E_i^\l f)(x)
\]
belongs to $\M_{\b-\a_i}$.
\end{proposition}

Denote by $E^\l_i$ the endomorphism of $\M$ mapping $f\in\M_\b$ to
$E_i^\l f$.
Notice that Formula~(\ref{actF}), which is nothing but
(\ref{stari}), also defines an endomorphism of $\M$ independent of
$\l$ which we again denote by $F_i$.
Finally Formula~(\ref{actH}) makes sense for any $\l$, 
not necessarily dominant, and any $f\in\M_\b$.
This gives an endomorphism of $\M$ that we shall denote by $H^\l_i$.

\begin{theorem}\label{conjV}
The assignments $e_i\mapsto E^\l_i$, $f_i\mapsto F_i$, $h_i\mapsto H^\l_i$,
give a representation of $\g$ on $\M$ isomorphic to the
Verma module $M(\l)$.
\end{theorem}

The rest of this section is devoted to the proofs of 
Proposition~\ref{conjf} and Theorem \ref{conjV}.

\subsection{}
Denote by $e^\l_i$ the endomorphism of the Verma module
$M(\l)$ implementing the action of the Chevalley generator $e_i$.
Let $\E^\l_i$ denote the endomorphism of $U(\nm)$
obtained by transporting $e^\l_i$ via the natural identification
$M(\l) \cong U(\nm)$.
Let $\Delta$ be the comultiplication of $U(\nm)$.
\begin{lemma}\label{lemalg1}
For $\l,\mu\in P$ and $u\in U(\nm)$ we have
\[
\Delta(\E^{\l+\mu}_i u) = 
(\E^\l_i\otimes 1 + 1\otimes\E^\mu_i)\Delta u .
\]
\end{lemma}
\proof
By linearity it is enough to prove this for $u$
of the form $u=f_{i_1}\cdots f_{i_r}$.
A simple calculation in $U(\g)$ shows that
\[
e_if_{i_1}\cdots f_{i_r} =
f_{i_1}\cdots f_{i_r}e_i +
\sum_{k=1}^r\de_{ii_k}f_{i_1}\cdots f_{i_{k-1}}h_if_{i_{k+1}}\cdots
f_{i_r} 
\]
\[
=f_{i_1}\cdots f_{i_r}e_i +
\sum_{k=1}^r\de_{ii_k}\left(f_{i_1}\cdots f_{i_{k-1}}f_{i_{k+1}}\cdots
f_{i_r}h_i -
\left(\sum_{s=k+1}^r a_{ii_s}\right)f_{i_1}\cdots f_{i_{k-1}}f_{i_{k+1}}\cdots
f_{i_r}\right).
\]
It follows that, for $\nu\in P$,
\[
\E^\nu_i(f_{i_1}\cdots f_{i_r}) =
\sum_{k=1}^r\de_{ii_k}\left((\nu;\a_i)-\sum_{s=k+1}^r a_{ii_s}\right)
f_{i_1}\cdots f_{i_{k-1}}f_{i_{k+1}}\cdots f_{i_r}.
\]
Now, using that $\De$ is the algebra homomorphism defined by
$\De(f_i)=f_i\otimes 1 + 1\otimes f_i$, one can finish the proof
of the lemma. Details are omitted.
\cqfd

\subsection{}
We endow $U(\nm)$ with the $Q_+$-grading given by 
$\deg(f_i)=\a_i$. 
Let $u$ be a homogeneous element of $U(\nm)$. 
Write $\De u = u \otimes 1 + u^{(i)} \otimes f_i + A$,
where $A$ is a sum of homogeneous terms of the form $u'\otimes u''$
with $\deg(u'') \not = \a_i$.
This defines $u^{(i)}$ unambiguously.
\begin{lemma}\label{lemalg2}
For $\l,\mu\in P$ we have
\[
\E^{\l+\mu}_i u = \E^\l_i u + (\mu;\a_i)\,u^{(i)}.
\]
\end{lemma}
\proof
We calculate in two ways the unique term of the form 
$E\otimes 1$ in $\Delta(\E^{\l+\mu}_i u)$.
On the one hand, we have obviously 
$E\otimes 1=\E^{\l+\mu}_i u \otimes 1$.
On the other hand, using Lemma~\ref{lemalg1}, we have
\[
E\otimes 1 = \E^{\l}_i u \otimes 1 
+ (1\otimes \E^{\mu}_i)(u^{(i)}\otimes f_i)
= \E^{\l}_i u \otimes 1 + (\mu;\a_i)\,u^{(i)}\otimes 1.
\]
Therefore,
\[
E=\E^{\l+\mu}_i u = \E^\l_i u + (\mu;\a_i)\,u^{(i)}.
\]
\cqfd

\subsection{}
Now let us return to the geometric realization $\M$ of $U(\nm)$.
Let $E^\l_i$ denote the endomorphism of $\M$
obtained by transporting $e^\l_i$ via the identification
$M(\l) \cong \M$.

\begin{lemma}\label{lem10}
Let $\l\in P_+$, $f\in\M_\b$ and $x\in\L^\l_{\b-\a_i}$.
Then
\[
(E^\l_if)(x) = \int_{y\in\G(x,\l,i)} f(y).
\]
\end{lemma}
\proof
Let $r_\l : \M \ra \M^\l$ be the linear map sending $f\in\M_\b$
to its restriction to $\L^\l_\b$.
By Theorem~\ref{thI}, this is a homomorphism of $U(\nm)$-modules
mapping the highest weight vector of $\M \cong M(\l)$ to the
highest weight vector of $\M^\l \cong L(\l)$.
It follows that $r_\l$ is in fact a homomorphism of $U(\g)$-modules,
hence the restriction of $E^\l_if$ to $\L^\l_{\b-\a_i}$
is given by Formula~(\ref{actE}) of Section~\ref{gri}.
\cqfd

Let again $\l\in P$ be arbitrary, and pick $f\in\M_\b$. 
It follows from Lemma~\ref{lemalg2}
that for any  $\mu\in P$
\[
E^{\l+\mu}_i f - (\mu;\a_i)\,f^{(i)}
=
E^\l_i f.
\]
Let $x\in\L_{\b-\a_i}$.
Choose $\nu=\l+\mu$ sufficiently dominant so that 
$x$ is isomorphic to a submodule of $q_\nu$. 
Then by Lemma~\ref{lem10}, we have
\[
(E^\nu_i f)(x) = \int_{y\in\G(x,\nu,i)} f(y).
\]
On the other hand, by the geometric description of $\De$
given in \cite[\S 6.1]{GLS}, if we write 
\[
\De f = f \otimes 1 + f^{(i)} \otimes \1_i + A
\]
where $A$ is a sum of homogeneous terms of the form $f'\otimes f''$
with $\deg(f'') \not = \a_i$, we have that $ f^{(i)}$ is the 
function on $\L_{\b-\a_i}$ given by 
$f^{(i)}(x) = f(x\oplus s_i)$.
Hence we obtain that for $x\in\L_{\b-\a_i}$ 
\[
(E^\l_i f)(x) = \int_{y\in\G(x,\nu,i)} f(y) \ - \ (\nu -
\l;\a_i)f(x\oplus s_i).
\]
This proves both Proposition~\ref{conjf} and Theorem~\ref{conjV}.
\cqfd

\subsection{}
Let $\l\in P_+$. We note the following consequence of
Lemma~\ref{lem10}.
\begin{proposition}
Let $\l\in P_+$. The linear map $r_\l : \M \ra \M^\l$ sending $f\in\M_\b$
to its restriction to $\L^\l_\b$ is the geometric realization of the
homomorphism of $\g$-modules $M(\l) \ra L(\l)$. \cqfd
\end{proposition}

\section{Dual Verma modules}

\subsection{}
Let $S$ be the anti-automorphism of $U(\g)$ defined by
\[
S(e_i) = f_i,\quad S(f_i)=e_i,\quad S(h_i) = h_i, \qquad (i\in I).
\]
Recall that, given a left $U(\g)$-module $M$, the dual module $M^*$
is defined by 
\[
(u\,\varphi)(m) = \varphi(S(u)\,m), \qquad 
(u\in U(\g),\ m\in M,\ \varphi\in M^*).
\]
This is also a left module.
If $M$ is an infinite-dimensional module with finite-dimensional
weight spaces $M_\nu$, we take for $M^*$ the graded dual
$M^*=\bigoplus_{\nu\in P} M_\nu^*$.

For $\l\in P$ we have $L(\l)^*\cong L(\l)$, hence the quotient
map $M(\l) \ra L(\l)$ gives by duality an embedding
$L(\l) \ra M(\l)^*$ of $U(\g)$-modules.

\subsection{}
Let $\M^* = \bigoplus_{\b\in Q_+} \M_\b^*$ denote the vector space 
graded dual of $\M$. 
For $x\in \L_\b$, we denote by $\de_x$ the delta function
given by 
\[
\de_x(f) = f(x),\qquad (f\in\M_\b).
\]
Note that the map $\de : x \mapsto \de_x$ is a constructible map from
$\L_\b$ to $\M_\b^*$. Indeed the preimage of $\de_x$ is the
intersection of the constructible subsets
\[
\M_{(i_1,\ldots,i_r)}=
\{y \in \L_\b \mid (\1_{i_1}*\cdots *\1_{i_r})(y)  
= (\1_{i_1}*\cdots *\1_{i_r})(x)\},
\quad (\a_{i_1}+\cdots+\a_{i_r}=\b).
\]

\subsection{}
We can now dualize the results of Sections~\ref{gri}
and \ref{grv} as follows.
For $\l\in P$ and $x\in\L_\b$ put
\begin{eqnarray}
(E_i^*)(\de_x) &=& \int_{y\in\G(i,x)} \de_y,
\label{actEs} \\[3mm]
(F_i^{\l*})(\de_x) &=& \int_{y\in\G(x,\nu,i)} \de_y
\ -\ (\nu-\l\,;\,\a_i)\,\de_{x\oplus s_i},
\label{actFs}\\[3mm]
(H_i^{\l*})(\de_x) &=& (\l-\b;\a_i)\,\de_x,
\label{actHs}
\end{eqnarray}
where in (\ref{actFs}) the weight $\nu\in P_+$ is such
that $x$ is isomorphic to a submodule of $q_\nu$.
The following theorem then follows immediately from Theorems~\ref{thI} and
\ref{conjV}.
\begin{theorem}\label{dual}
{\rm (i)}\ The formulas above define endomorphisms $E_i^*, F_i^{\l*}, H_i^{\l*}$
of $\M^*$, and the assignments $e_i\mapsto E^*_i$, $f_i\mapsto F^{\l*}_i$, 
$h_i\mapsto H^{\l*}_i$,
give a representation of $\g$ on $\M^*$ isomorphic to the
dual Verma module $M(\l)^*$.

{\rm (ii)}\ If $\l\in P_+$, the subspace $\M^{\l*}$ of $\M^*$ spanned
by the delta functions $\de_x$ of the finite-dimensional nilpotent
submodules $x$ of $q_\l$ carries the irreducible submodule
$L(\l)$. 
For such a module $x$, Formula~(\ref{actFs}) simplifies as follows
\[
(F_i^{\l*})(\de_x) = \int_{y\in\G(x,\l,i)} \de_y\,.
\]
\cqfd
\end{theorem}

\begin{example}
{\rm
Let $\g$ be of type $A_2$. 
Take $\l=\varpi_1+\varpi_2$, where
$\varpi_i$ is the fundamental weight corresponding to $i\in I$. 
Thus $L(\l)$ is isomorphic to the 8-dimensional adjoint 
representation of $\g=\Sl_3$.

A $\L$-module $x$ 
consists of a pair of linear maps $x_{21} : V_1 \ra V_2$ and 
$x_{12} : V_2 \ra V_1$ such that 
$x_{12}x_{21}=x_{21}x_{12}=0$.
The injective $\L$-module $q=q_\l$ has the following
form :
\[
q=
\begin{pmatrix}
u_1 \lra u_2 \cr v_1 \lla v_2
\end{pmatrix}
\]
This diagram means
that $(u_1,v_1)$ is a basis of $V_1$, that $(u_2,v_2)$ is a basis of 
$V_2$, and that 
\[
q_{21}(u_1)=u_2, \quad q_{21}(v_1)=0, \quad
q_{12}(v_2)=v_1, \quad q_{12}(u_2)=0.
\]
Using the same type of notation, we can exhibit the following submodules of $q$ :
\[
x_1=\begin{pmatrix}v_1\end{pmatrix}, 
\quad x_2=\begin{pmatrix}u_2\end{pmatrix}, \quad
x_3=\begin{pmatrix}v_1 && u_2\end{pmatrix}, \quad 
x_4=\begin{pmatrix}u_1 \lra u_2\end{pmatrix}, \quad 
x_5=\begin{pmatrix}v_1 \lla v_2\end{pmatrix}, 
\]
\[
x_6=\begin{pmatrix}u_1 \lra u_2 \cr v_1 \hfill\end{pmatrix},\qquad
x_7=\begin{pmatrix}\hfill u_2 \cr v_1 \lla v_2\end{pmatrix}.
\]
This is not an exhaustive list. For example, 
$x'_4=\begin{pmatrix}(u_1+v_1) \lra u_2\end{pmatrix}$ is another submodule, isomorphic 
to $x_4$.
Denoting by $\0$ the zero submodule, we see that $\de_\0$
is the highest weight vector of $L(\l)\subset M(\l)^*$.
Next, writing for simplicity $\de_i$ instead of $\de_{x_i}$
and $F_i$ instead of $F^\l_i$, Theorem~\ref{dual}~(ii) gives
the following formulas for the action of the $F_i$'s on 
$L(\l)$.
\[
F_1 \de_\0 = \de_1,\quad
F_2 \de_\0 = \de_2,\quad
F_1\de_2 = \de_3 + \de_4,\quad
F_2\de_1 = \de_3 + \de_5,
\]
\[
F_1 \de_3 = F_1 \de_4 =\de_6,\quad
F_2 \de_3 = F_2 \de_5 =\de_7,\quad
F_2 \de_3 = F_1 \de_6 = \de_q,\quad
F_1 \de_q = F_2 \de_q = 0.
\]
Now consider the $\L$-module $x = s_1 \oplus s_1$.
Since $x$ is not isomorphic to a submodule of $q_\l$,
the vector $\de_x$ does not belong to $L(\l)$.
Let us calculate $F_i \de_x\ (i = 1,2)$ by means of Formula~(\ref{actFs}).
We can take $\nu = 2\varpi_1$.
The injective $\L$-module $q_\nu$ has the following
form :
\[
q_\nu=
\begin{pmatrix}
w_1 \lla w_2 \cr v_1 \lla v_2
\end{pmatrix}
\]
It is easy to see that the variety $\G(x,\nu,2)$ is isomorphic
to a projective line $\P_1$, and that all points on this
line are isomorphic to 
\[
y=
\begin{pmatrix}
w_1\hfill \cr v_1 \lla v_2
\end{pmatrix}
\]
as $\L$-modules.
Hence,
\[
F_2 \de_x = \chi(\P_1)\,\de_y - (\nu - \l ; \a_2)\, \de_{x\oplus s_2}
= 2\,\de_y+\de_{s_1\oplus s_1\oplus s_2}.
\]
On the other hand, $\G(x,\nu,1)=\emptyset$, so that
\[
F_1 \de_x = - (\nu - \l ; \a_1)\, \de_{x\oplus s_1}
= - \de_{s_1\oplus s_1\oplus s_1}.
\]
}
\end{example}


\bigskip
\small

\noindent
\begin{tabular}{ll}
Christof {\sc Gei{ss}} : &
Instituto de Matem\'aticas, UNAM\\
& Ciudad Universitaria, 04510 Mexico D.F., Mexico \\
&email : {\tt christof@math.unam.mx}\\[5mm]
Bernard {\sc Leclerc} :&
LMNO,
Universit\'e de Caen,\\
& 14032 Caen cedex, France\\
&email : {\tt leclerc@math.unicaen.fr}\\[5mm]
Jan {\sc Schr\"oer} :&
Department of Pure Mathematics,
University of Leeds,\\
&Leeds LS2 9JT, England\\
&email : {\tt jschroer@maths.leeds.ac.uk}
\end{tabular}


\begin{thebibliography}{ABC} \scriptsize
%
\bibitem[\bf Bo]{Bo}
{\sc K. Bongartz},
{\it Algebras and quadratic forms}, 
J. London Math. Soc. {\bf 28} (1983), 461--469. 

\bibitem[\bf BK]{BK}
{\sc M. C. R. Butler, A. D. King}, 
{\it Minimal resolutions of algebras}, 
J. Algebra {\bf 212} (1999), 323--362. 
%
\bibitem[\bf CB]{CB}
{\sc W. Crawley-Boevey}, 
{\it On the exceptional fibres of Kleinian singularities},
Amer. J. Math. {\bf 122} (2000), 1027--1037.
%
\bibitem[\bf GLS]{GLS}{\sc C. Geiss, B. Leclerc, J. Schr\"oer},
{\it Semicanonical bases and preprojective algebras},
Ann. Scient. \'Ec. Norm. Sup. {\bf 38} (2005), 193--253.
%
\bibitem[\bf GP]{GP}{\sc I. M. Gelfand, V. A. Ponomarev},
{\it Model algebras and representations of graphs},
Funct. Anal. Appl. {\bf 13} (1980), 157--166.
%
\bibitem[\bf Lu1]{Lu91}{\sc G. Lusztig},
{\it Quivers, perverse sheaves, and quantized enveloping
algebras}, J. Amer. Math. Soc. {\bf 4} (1991), 365--421. 
%
\bibitem[\bf Lu2]{Lu00}{\sc G. Lusztig}, 
{\it Semicanonical bases arising from enveloping algebras},
Adv. Math. {\bf 151} (2000), 129--139.
%
\bibitem[\bf Lu3]{Lu002}{\sc G. Lusztig}, 
{\it Remarks on quiver varieties},
Duke Math. J. {\bf 105} (2000), 239--265.

\bibitem[\bf Lu4]{Lu03}{\sc G. Lusztig}, 
{\it Constructible functions on varieties attached to quivers},
in {Studies in memory of Issai Schur} 177--223, 
Progress in Mathematics {\bf 210}, Birkh\"auser 2003.
%
\bibitem[\bf Na]{Na}{\sc H. Nakajima},
{\it Instantons on ALE spaces, quiver varieties, and
Kac-Moody algebras},
Duke Math. J. {\bf 76} (1994), 365--416.
%
\bibitem[\bf Ri]{Ri}
{\sc C. M. Ringel},
{\it The preprojective algebra of a quiver}, 
in {\it Algebras and modules} II (Geiranger, 1966),
467--480, CMS Conf. Proc. {\bf 24}, AMS 1998.
%
\end{thebibliography}
\end{document}